\documentclass[12pt]{article}

\usepackage{geometry} 
\usepackage{amsthm,amssymb}
\usepackage{amsmath,amsfonts,latexsym}
\usepackage{latexsym}
\usepackage{float}
\usepackage{xcolor}
\usepackage{url}

\def\titlerunning#1{\gdef\titrun{#1}}
\makeatletter
\def\author#1{\gdef\autrun{\def\and{\unskip, }#1}\gdef\@author{#1}}
\def\address#1{{\def\and{\\\hspace*{18pt}}\renewcommand{\thefootnote}{}%
\footnote {#1}}%
\markboth{\autrun}{\titrun}}
\makeatother
\def\email#1{\hspace*{4pt}{\em e-mail}: #1}

\newtheorem{thm}{Theorem}[section]

\theoremstyle{definition}

\begin{document}

\titlerunning{}

\title{Strongly regular graphs with parameters $(81,30,9,12)$ and a new partial geometry}

\author{Dean Crnkovi\' c, Andrea \v{S}vob and Vladimir D. Tonchev}

\maketitle

\address{D. Crnkovi\' c, A. \v{S}vob: Department of Mathematics, University of Rijeka, Radmile Matej\v{c}i\'c 2, 51000 Rijeka, Croatia;
\email{\{deanc,asvob\}@math.uniri.hr}
\and
V. D. Tonchev: Department of Mathematical Sciences, Michigan Technological University, Houghton, MI 49931, USA;
\email{tonchev@mtu.edu} 
}


\begin{abstract}
Twelve new strongly regular graphs with parameters 
$(81,30,9,12)$ are found as graphs invariant under
certain subgroups of the automorphism groups of
the two previously known graphs that arise from 2-weight codes.
One of these new graphs is geometric and yields a partial geometry  with parameters $pg(5,5,2)$ that is not isomorphic to the partial geometry  discovered by J. H. van Lint and A. Schrijver 
\cite{van-Lint} in 1981.
\end{abstract} 

\section{Introduction}

A {\em strongly regular graph} with parameters $(v, k, \lambda, \mu)$ (or $srg(v, k, \lambda, \mu)$ for short) is an undirected  graph with $v$ verices, having no multiple edges or loops,
such that every vertex has exactly  $k$ neighbours, every two adjacent vertices have exactly 
$\lambda$ common neighbours, and every two non--adjacent vertices have exactly $\mu$ common neighbours.

A {\it partial geometry} with parameters $s, t, \alpha$,
or shorty,  $pg(s,t,\alpha)$, is a pair $(P,L)$ of a set $P$
of {\it points} and a set $L$ of {\it lines}, with an incidence relation
between points and lines, satisfying the following axioms:
\begin{enumerate}
\item A pair of distinct points is not incident with more than one
line.
\item Every line  is incident with exactly $s+1$ points ($s\ge 1$).
\item Every point is incident with exactly $t+1$ lines ($t\ge 1$).
\item For every point $p$ not incident with a line $l$, 
there are exactly $\alpha$ lines ($\alpha \ge 1$) 
which are incident with $p$, 
and also incident with some point incident with $l$.
\end{enumerate} 

Partial geometries and strongly regular graphs were introduced by R. C. Bose \cite{Bose}.
In the original Bose's notation, the 
number $t+1$ of lines incident with a point 
is denoted by $r$, and the number $s+1$ of points incident 
with a line is denoted by $k$. 
A survey on strongly regular graphs is given by Brouwer in 
\cite{crc-srg},
and for a survey on partial geometries, see Thas \cite{crc-pg}.

In terms of $s, t, \alpha$,
the  number $v=|P|$ of  points, and the number
$b=|L|$ of lines of a partial geometry $pg(s,t,\alpha)$
are given by eq. (\ref{eqvb}).
\begin{equation}
\label{eqvb}
v=\frac{(s+1)(st + \alpha)}{\alpha}, \ b=\frac{(t+1)(st + \alpha)}{\alpha}. 
\end{equation}
If $G=(P,L)$ is a partial geometry $pg(s,t,\alpha)$,
the incidence structure $G'$ having as points the lines of $G$, 
and having as lines the points of $G$, where a point and a line
are incident in $G'$ if and only if the corresponding line and a point 
of $G$ are incident, is a partial geometry $pg(t,s,\alpha)$, called the {\it dual} of $G$. 

If $G=(P,L)$ is a partial geometry $pg(s,t,\alpha)$
 with point set $P$ and line set $L$,
the {\it point graph} $\Gamma_P$ of $G$ is the graph with vertex set $P$,
where two vertices are adjacent if the corresponding points of $G$ are
collinear. The {\it line graph} $\Gamma_L$ of $G$
is the graph having as vertices the lines of $G$, where two lines are adjacent
if they share a point. Both $\Gamma_P$ and $\Gamma_L$ are strongly regular 
graphs \cite{Bose}. The parameters $v, k, \lambda, \mu$ of $\Gamma_P$
are given by (\ref{eqP}).
\begin{equation}
\label{eqP}
v=(s+1)(st+\alpha)/\alpha,  k = s(t+1),  \lambda = s-1 +t(\alpha -1), \mu = \alpha(t+1). 
\end{equation}
A strongly regular graph $\Gamma$ whose parameters
$n, k, \lambda, \mu$ can be written as in eq. (\ref{eqP}) for some
integers $s, t, \alpha$ is called {\it pseudo-geometric},
and  $\Gamma$ is called {\it geometric}
 if there exists a partial geometry $G$
with parameters $s, t, \alpha$ such that $\Gamma$ is the point graph  of $G$.
A pseudo-geometric strongly regular graph $\Gamma$ with parameters (\ref{eqP}) is geometric if and only if there exists a set $S$ of $b=(t+1)(st + \alpha)/\alpha$ cliques of size $s+1$ such that every two cliques from $S$ share at most one vertex
\cite{Bose}.

A partial geometry $pg(s,t,\alpha)$ with $\alpha=s+1$
 is a Steiner 2-$(v,s+1,1)$ design, or dually, if $\alpha=t+1$,
 then the dual geometry is a Steiner 2-$(b,t+1,1)$ design.
 If $\alpha=s$, or dually, $\alpha=t$, then $G$ is a net of order 
 $s+1$ and degree $t+1$ \cite{Bose}. 
 A partial geometry with $\alpha = 1$ is a generalized quadrangle
\cite{Bose},  \cite{PT}.
 
 A partial geometry $pg(s,t,\alpha)$ is called {\it proper}
 if $1< \alpha < min(s,t)$.
 
The known proper partial geometries are divided into eight 
types,
four of which are infinite families, and there are four sporadic geometries that do not belong to any known  infinite family
\cite[Theorem 41.31]{crc-pg}.
One of the four sporadic examples is a partial geometry 
$pg(5,5,2)$   discovered by van Lint and Schrijver \cite{van-Lint}.
The point graph of the van Lint-Schrijver geometry has
parameters $v=81, \ k=30, \lambda=9, \mu =12$,
and is invariant under the elementary abelian group of order 81
acting regularly on the set of vertices.
By a result of Delsarte \cite{Del}, any graph with these parameters
that is  invariant under the elementary abelian group of order 81
acting regularly on the set of vertices,
can be obtained from 
a ternary linear $[30,4,9]$  two-weight code
with weight distribution $a_9 =50, a_{12}=30$
(see also \cite{aeb}, \cite{two-weight}, \cite{CalKan}).
Up to equivalence, there are exactly two such codes
(Hamada and Helleseth \cite{HamHel}), that give rise
to two nonisomorphic $srg(81,30,12,9)$:
a graph $\Gamma_1$,  being
isomorphic to the  van Lint-Schrijver graph and having full 
automorphism group of order    116640 ,
and a  second graph $\Gamma_2$  having full automorphism group of order 5832.
According to \cite{two-weight}, \cite{aeb},
\cite{aebp},  $\Gamma_1$ and $\Gamma_2$ appear to be the only previously known strongly regular graphs with parameters 
$v=81, k=30, \lambda=9, \mu=12$  .

In this paper, we use a method for finding strongly regular graphs based on orbit matricest that was developed in \cite{BL}   and \cite{CM}, to  show that  there are exactly three nonisomorphic graphs
$srg(81,30,9,12)$ which are invariant under a subgroup of order 360 of the automorphism group of $\Gamma_1$,  and exactly  eleven nonisomorphic graphs  $srg(81,30,9,12)$  which are invariant under a subgroup of order 972 of the automorphism group of $\Gamma_2$. One of the newly found graphs invariant under a group of order 972 gives rise to a new partial geometry $pg(5,5,2)$ that is not isomorphic to the van Lint-Schrijver  partial geometry. 
An isomorphic partial geometry was simultaneously and independently constructed by V. Kr\v cadinac  \cite{vk-preprint} by using a different method.
The adjacency matrices of the two previously known graphs and the twelve newly found graphs are available online at
\begin{verbatim}
 http://www.math.uniri.hr/~asvob/SRGs81_pg552.txt
\end{verbatim}
The  lines of the new  are given in the Appendix.

\section{New strongly regular graphs  $srg(81,30,9,12)$}\label{srgs}

 The graphs $\Gamma_1$ and $\Gamma_2$ described in
 the introduction
  have full automorphism groups $G_1$ and $G_2$ of order 
116640 and 5832, respectively. 
In this section, we construct twelve new  strongly regular graphs
 with parameters 
$(81,30,9,12)$. These new graphs  are constructed by 
expanding orbit matrices with respect to the action of certain subgroups of $G_1$ or $G_2$.
For more information on orbit matrices of  strongly regular graphs
 we refer the reader to \cite{BL}, \cite{CM}.
We used Magma for all computations involving groups and codes  in this paper.
\cite{magma}.

\subsection{Graphs invariant under  subgroups $A_6$ of 
$\mathrm{Aut}(\Gamma_1)$} \label{srg116640}

There are exactly four conjugacy classes of subgroups of order 360 in the group $G_1=\mathrm{Aut}(\Gamma_1)$ of order 116640, the representatives of which will be denoted by $H_1^1,\dots, H_4^1$. 
Each of these four representatives is isomorphic to the simple group $A_6$.

The subgroup $H_1^1$ is acting in two orbits on the set of vertices of $\Gamma_1$, one of size 36 and the other of size 45, 
giving an  orbit matrix $OM_1^1$ (\ref{om11}).
\begin{equation}  
\label{om11}
OM_1^1=   \left(
\begin{tabular}{rr}
15  & 15 \\
12  & 18 \\
\end{tabular}  \right)
\end{equation}

The orbit matrix $OM_1^1$  
expands to the  (0,1)-adjacency matrices
of two 
non-isomorphic strongly regular graphs:
 the graph $\Gamma_{1}$, and a new graph,
  denoted by $\Gamma_{14}$ with  full automorphism 
  group of order 360. 

The subgroup $H_2^1$ is acting in three orbits of sizes 6, 15 
and 60
respectively, giving an orbit matrix $OM_2^1$.

\begin{displaymath}  OM_2^1=   \left(
\begin{tabular}{rrr}
0  & 0  &  30 \\
0  & 6  &  24 \\
3  & 6  &  21 \\
\end{tabular} \right)
\end{displaymath}

The orbit matrix $OM_2^1$ gives rise  to  two non-isomorphic strongly regular graphs, $\Gamma_{1}$, and  a second graph  
 $\Gamma_{13}$ having  full automorphism group of order 720. 

The subgroup $H_3^1$ is acting in three orbits with sizes 6, 15 and 60, and  orbit matrix $OM_3^1$. 

\begin{displaymath}  OM_3^1=   \left(
\begin{tabular}{rrr}
5  & 5  &  20 \\
2  & 8  &  20 \\
2  & 5  &  23 \\
\end{tabular} \right)
\end{displaymath}

The orbit matrix  $OM_3^1$ can be expanded to only one
(up to isomorphism) strongly regular graph,   namely, the 
original graph 
$\Gamma_1$. 

The subgroup $H_4^1$  acts in four orbits, with sizes 1, 20, 30 and 30, and  orbit matrix $OM_4^1$. 

\begin{displaymath}  OM_4^1=   \left(
\begin{tabular}{rrrr}
0  & 0  &  0 & 30 \\
0  & 9  &  9 & 12 \\
0  & 6  &  12  &12 \\
1  & 8  &  12  &9 \\
\end{tabular} \right)
\end{displaymath}

Up to isomorphism, the matrix $OM_4^1$  is the orbit matrix
of only one strongly regular graph, that is, $\Gamma_1$. 

\subsection{Graphs invariant under subgroups of order 972 } \label{srg5832}

There are exactly five conjugacy classes of subgroups of order 972 in the group $G_2=\mathrm{Aut}(\Gamma_2)$ of order 5832, with representatives 
$H_1^2,\dots,H_5^2$.

The subgroups $H_1^2$ and $H_2^2$ act transitively on the 81 vertices and produce two non-isomorphic strongly regular graphs, isomorphic to $\Gamma_1$ and $\Gamma_2$. 

The subgroup $H_3^2$   partitions the set of vertices of 
$\Gamma_2$ in three orbits of length 27, and gives
an orbit matrix $OM_3^2$. 

\begin{displaymath}  OM_3^2=   \left(
\begin{tabular}{rrr}
12  & 9  &  9 \\
9  & 12  &  9 \\
9  & 9  &  12 \\
\end{tabular} \right)
\end{displaymath}

The orbit matrix $OM_3^2$ gives rise to two non-isomorphic strongly regular graphs, the graph $\Gamma_2$, and a new graph $\Gamma_4$, having  full automorphism group of order 1944.

The subgroups $H_4^2$ and $H_5^2$ are acting in the same way, with two orbits of size 27 and  54 respectively,
 and gives an  orbit matrix $OM_4^2$. 
 
 \begin{displaymath}  
OM_4^2=   \left(
\begin{tabular}{rr}
12  & 18 \\
9  & 21 \\
\end{tabular}  \right)
\end{displaymath}

The subgroup $H_4^2$ leads to four non-isomorphic strongly regular graphs, incuding $\Gamma_2$. Two of these 
graphs, $\Gamma_5$ and $\Gamma_6$, have full automorphism groups of order 972, 
and $\Gamma_4$ has full automorphism group of order 1944.
The group $H_5^2$ gives rise to 12 non-isomorphic strongly regular graphs, among them $\Gamma_1$, $\Gamma_2$ and 
$\Gamma_4$. 
Eight of these graphs, denoted by $\Gamma_5,\dots, \Gamma_{12}$, have  full automorphism groups of order 972, 
and $\Gamma_3$ has  full automorphism group of order 3888. 

The graph $\Gamma_{12}$ is geometric and produce a new partial geometry $pg(5,5,2)$.

\section{A new partial geometry $pg(5,5,2)$}

The results presented in Section \ref{srgs} can be summarized as follows.

\begin{thm}
\begin{enumerate}
\item Up to isomorphism, there are exactly $3$ strongly regular graphs with parameters $(81,30,9,12)$ invariant under a subgroup
of order $360$ of the automorhism group of the graph 
$\Gamma_1$.

\item Up to isomorphism, there are exactly twelve strongly regular graphs with parameters $(81,30,9,12)$ invariant under a subgroup of
order $972$ of the automorphism group of the graph 
$\Gamma_2$.
\item One of the twelve graphs from part (2) yeilds a new
partial geometry $pg(5,5,2)$.
\end{enumerate}

\end{thm}
 
Details about these strongly regular graphs are given in 
Table \ref{results}:
the order  of the full automorphism group
of the graph, the 3-rank of the $(0,1)$-adjacency matrix,
the maximum clique size, and  the number of 6-cliques.
For every graph $\Gamma_i$ that contains at least 81 6-cliques,
 we define a graph $\Gamma_{i}^*$
having as vertices the 6-cliques of $\Gamma_i$, where two 6-cliques
are adjacent if there share at most one vertex.
The last but one column of Table \ref{results} contains the maximum
clique size of $\Gamma_{i}^*$, and if this maximum clique size is 81,
the last column contains the total number of 81-cliques in  
$\Gamma_{i}^*$.

Only two of the fourteen graphs, $\Gamma_1$ and $\Gamma_{12}$,
 are geometric.
The graph $\Gamma_{1}^*$ contains two 81-cliques, each consisting
of 81 6-cliques of $\Gamma_1$ that are the lines of a partial
geometry $pg(5,5,2)$ with full automorphism
group of order $58320$ acting transitively on the 
sets of points and lines, 
and is isomorphic to the van Lint-Schrijver partial geometry 
\cite{van-Lint}.

The graph $\Gamma_{12}^*$ contains only one 81-clique
and  yields a new partial geometry $pg^*(5,5,2)$ 
that is not isomorphic to the van Lint-Schrijver partial geometry.
The full automorphism group of the new partial geometry
is of order 972, and partitions the set of points, as well
as the set of lines, in two orbits of length 54 and 27.
An automorphism $f\in Aut(\Gamma_{12})$ of order 6 and a 
set of orbit representatives of the lines of $pg^*(5,5,2)$ under the action
of $<f>$ are listed in Table \ref{reps}, where the first
three orbits are of length 3, and the remaining 12 orbits are
of length 6.
The set of all 81 lines of $pg^*(5,5,2)$  is given in the Appendix.

\begin{table}[H]
\begin{center} \begin{footnotesize}
\begin{tabular}{|c|}
\hline
A generator $f$ of order 6\\
\hline
(1, 4, 3)(2, 6, 9, 18, 7, 13)(5, 20, 10, 21, 11, 12)(8, 24, 19, 16, 15, 27)\\
(14, 17, 23, 26, 22, 25)(28, 29, 32)(30, 47, 42, 39, 36, 34)(31, 33, 35, 46, 37, 38)\\
(40, 41, 45, 54, 49, 50)(43, 48, 51, 53, 44, 52)(55, 71, 72, 65, 59, 57)\\
(56, 67, 77, 68, 63, 73)(58, 60, 75, 74, 62, 61)(64, 78, 66)(69, 81, 76, 79, 70, 80)\\
\hline
\hline
Line orbits representatives\\
\hline
\{ 4, 15, 24, 28, 31, 46 \}, \{ 32, 51, 52, 64, 70, 81 \}, \{ 4, 17, 22, 57, 72, 78 \},\\
\{ 31, 39, 41, 58, 67, 71 \}, \{ 2, 11, 24, 66, 71, 72 \}, \{ 7, 14, 20, 36, 43, 54 \},\\
\{ 5, 13, 26, 36, 44, 50 \}, \{ 29, 43, 52, 62, 68, 71 \}, \{ 9, 11, 27, 68, 74, 80 \},\\
\{ 33, 42, 49, 71, 73, 75 \}, \{ 35, 41, 47, 66, 76, 81 \}, \{ 1, 19, 24, 30, 44, 54 \},\\
\{ 9, 10, 24, 56, 62, 70 \}, \{ 4, 14, 25, 56, 69, 75 \}, \{ 11, 13, 25, 32, 37, 46 \}\\
\hline
\end{tabular}\end{footnotesize} 
\caption{\footnotesize The new  partial geometry $pg^*(5,5,2)$}\label{reps}
\end{center}
\end{table}

The data in Table \ref{results} distinguishes as nonisomorphic
all but  the three graphs $\Gamma_5, \Gamma_7, \Gamma_8$,
and the two graphs $\Gamma_9$ and $\Gamma_{10}$.
Let $D_i$ be the design on 81 points having as blocks
the 6-cliques in $\Gamma_i$. 
We checked with Magma \cite{magma} that $|Aut(D_9)|=972$, while 
$|Aut(D_{10}|=1944$,
hence $\Gamma_9$ and $\Gamma_{10}$ are nonisomorphic.
Similarly, $|Aut(D_5)|=972$, $|Aut(D_7)|=|Aut(D_8)|=1944$
shows that $\Gamma_5$ is not isomorphic to $\Gamma_7$ or 
$\Gamma_8$.
Finally, $\Gamma_7$ and $\Gamma_8$ can be shown to be nonisomorphic by comparing the weight distributions of the ternary linear codes spanned by their adjacency matrices.
The ternary code of $\Gamma_7$ contains 32400 codewords of
weight 33, while the code of $\Gamma_8$  contains 44550
codewords of weight 33.

\begin{table}[H]
\begin{center} \begin{footnotesize}
\begin{tabular}{|c|c|c|c|c|c|c|}
\hline
Graph $\Gamma$ &  $|Aut(\Gamma)|$& 3-rank & Max. clique size of $\Gamma$ & \# 6-cliques & Max. clique size of $\Gamma^*$ & \# 81-cliques\\
\hline
\hline
$\Gamma_1$ & 116640 & 19  &6 &162 &81 & 2 \\
$\Gamma_2$ & 5832 &19    &4 &0   &  & \\
$\Gamma_3$ & 3888&21   &6 &54  &  & \\
$\Gamma_4$ & 1944 &21    &6 &108 &54 &\\
$\Gamma_5$ & 972 &21     &6 &54  &  & \\
$\Gamma_6$ & 972 &21    &6 &108 &54 &\\
$\Gamma_7$ & 972 &21    &6 &54  &  & \\
$\Gamma_8$ & 972 &21    &6 &54  &  & \\
$\Gamma_9$ & 972&20     &6 &81  &54 &\\
$\Gamma_{10}$ & 972 &20  &6 &81  &54 & \\
$\Gamma_{11}$ & 972&21   &6 &108 &54 &\\
$\Gamma_{12}$ & 972 &21 &6 &108 &81 & 1 \\
$\Gamma_{13}$ & 720 &21 &6 &90  &45 &\\
$\Gamma_{14}$ & 360&25  &6 &21  &   &\\
\hline
\hline
\end{tabular}\end{footnotesize} 
\caption{\footnotesize Strongly regular graphs with parameters
 $(81,30,9,12)$}\label{results}
\end{center}
\end{table}

\noindent {\bf Acknowledgement} \\
D. Crnkovi\' c and A. \v Svob were supported by {\rm C}roatian Science Foundation under the project 6732.

\newpage
\section*{Appendix}

{\sc Lines of the new partial geometry $pg^*(5,5,2)$}

\begin{align*}
&\{\{ 1, 8, 16, 32, 35, 38 \}, \{ 1, 14, 26, 59, 64, 71 \}, \{ 1, 15, 27, 39, 40, 48 \}, \{ 1, 17, 23, 62, 76, 77 \}, \\
&\{ 1, 19, 24, 30, 44, 54 \}, \{ 1, 22, 25, 60, 73, 80 \}, \{ 2, 5, 27, 63, 75, 76 \}, \{ 2, 10, 16, 60, 67, 79 \}, \\
&\{ 2, 11, 24, 66, 71, 72 \}, \{ 2, 12, 14, 28, 37, 38 \}, \{ 2, 20, 22, 34, 40, 52 \}, \{ 2, 21, 23, 30, 50, 51 \}, \\
&\{ 3, 8, 24, 34, 45, 53 \}, \{ 3, 14, 17, 67, 74, 81 \}, \{ 3, 15, 16, 42, 43, 50 \}, \{ 3, 19, 27, 29, 33, 37 \}, \\
&\{ 3, 22, 26, 58, 63, 70 \}, \{ 3, 23, 25, 55, 65, 66 \}, \{ 4, 8, 27, 36, 41, 51 \}, \{ 4, 14, 25, 56, 69, 75 \}, \\
&\{ 4, 15, 24, 28, 31, 46 \}, \{ 4, 16, 19, 47, 49, 52 \}, \{ 4, 17, 22, 57, 72, 78 \},  \{ 4, 23, 26, 61, 68, 79 \},\\
 &\{ 5, 6, 17, 29, 31, 38 \}, \{ 5, 7, 24, 61, 73, 81 \}, \{ 5, 9, 16, 59, 65, 78 \}, \{ 5, 13, 26, 36, 44, 50 \}, \\
&\{ 5, 18, 25, 39, 45, 52 \}, \{ 6, 8, 20, 73, 74, 79 \}, \{ 6, 10, 25, 30, 41, 43 \}, \{ 6, 11, 26, 40, 47, 53 \}, \\
&\{ 6, 12, 19, 64, 65, 72 \}, \{ 6, 15, 21, 70, 75, 77 \}, \{ 7, 10, 27, 55, 57, 64 \}, \{ 7, 11, 16, 58, 69, 77 \}, \\
&\{ 7, 12, 23, 39, 49, 53 \}, \{ 7, 14, 20, 36, 43, 54 \}, \{ 7, 21, 22, 29, 35, 46 \}, \{ 8, 12, 18, 62, 63, 69 \},\\
 &\{ 8, 13, 21, 55, 71, 78 \}, \{ 9, 10, 24, 56, 62, 70 \}, \{ 9, 11, 27, 68, 74, 80 \}, \{ 9, 12, 22, 41, 42, 44 \}, \\
&\{ 9, 14, 21, 45, 47, 48 \}, \{ 9, 20, 23, 31, 32, 33 \}, \{ 10, 13, 17, 34, 48, 49 \}, \{ 10, 18, 26, 28, 33, 35 \}, \\
&\{ 11, 13, 25, 32, 37, 46 \}, \{ 11, 17, 18, 42, 51, 54 \}, \{ 12, 13, 15, 60, 68, 81 \}, \{ 13, 19, 20, 56, 58, 76 \}, \\
&\{ 15, 18, 20, 57, 59, 66 \},  \{ 18, 19, 21, 61, 67, 80 \}, \{ 28, 43, 53, 76, 78, 80 \}, \{ 28, 44, 52, 55, 74, 77 \},\\
& \{ 28, 48, 51, 58, 65, 73 \}, \{ 29, 43, 52, 62, 68, 71 \}, \{ 29, 44, 48, 66, 69, 79 \}, \{ 29, 51, 53, 56, 59, 60 \}, \\
&\{ 30, 33, 40, 69, 78, 81 \}, \{ 30, 38, 45, 57, 58, 68 \}, \{ 30, 46, 49, 59, 63, 74 \}, \{ 31, 34, 50, 64, 69, 80 \}, \\
&\{ 31, 39, 41, 58, 67, 71 \}, \{ 31, 47, 54, 55, 60, 63 \}, \{ 32, 43, 48, 61, 63, 72 \}, \{ 32, 44, 53, 57, 67, 75 \}, \\
&\{ 32, 51, 52, 64, 70, 81 \}, \{ 33, 36, 45, 60, 72, 77 \}, \{ 33, 42, 49, 71, 73, 75 \}, \{ 34, 35, 54, 65, 68, 75 \}, \\
&\{ 34, 37, 41, 59, 61, 77 \}, \{ 35, 39, 50, 56, 72, 74 \}, \{ 35, 41, 47, 66, 76, 81 \}, \{ 36, 38, 49, 66, 70, 80 \}, \\
&\{ 36, 40, 46, 62, 65, 67 \}, \{ 37, 39, 54, 70, 78, 79 \}, \{ 37, 47, 50, 57, 62, 73 \}, \{ 38, 40, 42, 55, 56, 61 \}, \\
&\{ 42, 45, 46, 64, 76, 79 \} \}
\end{align*}

\end{document}